\documentclass[11pt,leqno,twoside]{article}

\usepackage{amsfonts,amsmath,amsthm,amssymb}
\usepackage{enumerate}
\usepackage{graphics,graphicx,subfigure}

\usepackage{color}
\usepackage{todonotes}
\usepackage{cancel}
\usepackage{url}
\usepackage{hyperref}
\usepackage{makeidx}
\usepackage{showidx}
\usepackage{multicol}        
\usepackage{xspace}
\usepackage{stmaryrd}        
\usepackage{pifont}          
\usepackage{fancybox}        
\usepackage{bm}

\usepackage{fancyhdr}
\theoremstyle{plain}
\theoremstyle{definition}
\newtheorem{lem}{Lemma}
\newtheorem{defn}[lem]{Definition}
\newtheorem{thm}[lem]{Theorem}
\newtheorem{prop}[lem]{Proposition}
\newtheorem{cor}[lem]{Corollary}
\newtheorem{notn}[lem]{Notations}
\newtheorem{pb}[lem]{Problem}
\newtheorem{form}[lem]{Formulae}

\newtheorem*{rk}{Remark}
\newtheorem*{com}{Comment}
\newtheorem*{ex}{Example}
\theoremstyle{remark}

\newcommand{\blem}{\begin{lem}}
	\newcommand{\elem}{\end{lem}}
\newcommand{\bdefn}{\begin{defn}}
	\newcommand{\edefn}{\end{defn}}
\newcommand{\bthm}{\begin{thm} }
	\newcommand{\ethm}{\end{thm}}
\newcommand{\bprop}{\begin{prop}}
	\newcommand{\eprop}{\end{prop}}
\newcommand{\bcor}{\begin{cor}}
	\newcommand{\ecor}{\end{cor}}
\newcommand{\bnotn}{\begin{notn}}
	\newcommand{\enotn}{\end{notn}}
\newcommand{\bpb}{\begin{pb}}
	\newcommand{\epb}{\end{pb}}
\newcommand{\bform}{\begin{form}}
	\newcommand{\eform}{\end{form}}
\newcommand{\brk}{\begin{rk}}
	\newcommand{\erk}{\end{rk}}
\newcommand{\bcom}{\begin{com}}
	\newcommand{\ecom}{\end{com}}
\newcommand{\bex}{\begin{ex}}
	\newcommand{\eex}{\end{ex}}
\newcommand{\bpf}{\begin{proof}}
	\newcommand{\epf}{\end{proof}}

\newcommand{\va}{{\bf a}}
\newcommand{\vb}{{\bf b}}

\newcommand{\ve}{{\bf e}}

\newcommand{\vu}{{\bf u}}
\newcommand{\vv}{{\bf v}}
\newcommand{\vw}{{\bf w}}

\newcommand{\vy}{{\bf y}}
\newcommand{\vz}{{\bf z}}
\newcommand{\vA}{{\bf A}}
\newcommand{\vB}{{\bf B}}

\newcommand{\vX}{{\bf X}}

\newcommand{\vZ}{{\bf Z}}

\newcommand{\cA}{\mathcal{A}}

\newcommand{\cG}{\mathcal{G}}

\newcommand{\cN}{\mathcal{N}}

\newcommand{\cT}{\mathcal{T}}

\newcommand{\bR}{\mathbb{R}}

\newcommand{\be}{\begin{equation}}
	\newcommand{\ee}{\end{equation}}
\newcommand{\bal}{\begin{align}}
	\newcommand{\eal}{\end{align}}
\newcommand{\ba}{\begin{align*}}
	\newcommand{\ea}{\end{align*}}
\newcommand{\bmx}{\begin{matrix}}
	\newcommand{\emx}{\end{matrix}}
\newcommand{\bbmx}{\begin{bmatrix}}
	\newcommand{\ebmx}{\end{bmatrix}}
\newcommand{\bpmx}{\begin{pmatrix}}
	\newcommand{\epmx}{\end{pmatrix}}
\newcommand{\bvmx}{\begin{vmatrix}}
	\newcommand{\evmx}{\end{vmatrix}}

\newcommand{\f}{\frac}
\newcommand{\df}{\dfrac}

\newcommand{\tto}{\longrightarrow}

\newcommand{\tr}{\mathrm{tr}}
\newcommand{\sgn}{\mathrm{sgn}}

\newcommand{\minimize}[1]{\underset{#1}{\rm minimize}\,}

\newcommand{\la}{\lambda}

\newcommand{\eps}{\varepsilon}

\usepackage{mathtools}
\newcommand{\rev}[1]{#1} 
\newcommand{\revv}[1]{#1} 

\title{\vspace{-25mm}Iterative Hard Thresholding for Low-Rank Recovery\\ from Rank-One Projections\medskip\hrule height 1.2pt \vspace{-6mm}}
\author{Simon Foucart\footnote{S. F. partially supported by NSF grants DMS-1622134 and DMS-1664803.}  \, and Srinivas Subramanian --- Texas A\&M University}
\date{\vspace{-6mm}\rule{100mm}{0.8pt}}

\newcommand\shorttitle{Iterative Hard Thresholding for Low-Rank Recovery from Rank-One Projections}
\newcommand\authors{S. Foucart, S. Subramanian}

\begin{document}
\maketitle

\vspace{-15mm}
\begin{abstract}
A novel algorithm for the recovery of low-rank matrices acquired via compressive linear measurements is proposed and analyzed.
The algorithm,
a variation on the iterative hard thresholding algorithm for low-rank recovery,
is designed to succeed in situations where the standard rank-restricted isometry property fails,
e.g. in case of subexponential unstructured measurements
or of subgaussian rank-one measurements.
The stability and robustness of the algorithm are established based on distinctive matrix-analytic ingredients
and its performance is substantiated numerically.
\end{abstract}



\vspace{-5mm}
\begin{center}
\rule{100mm}{0.8pt}
\end{center}


\section{Introduction}

This article is concerned with the recovery of matrices $\vX \in \bR^{N_1 \times N_2}$ 
of rank $r \ll N:=\min\{N_1,N_2\}$ from compressive linear measurements
\be
\vy = \cA(\vX) \in \bR^{m}
\ee
with a number $m$ of measurements scaling like $r \max\{N_1,N_2\}$ instead of  $N_1N_2$.
It is known that such a low-rank recovery task can be successfully accomplished 
when the linear map $\cA: \bR^{N_1 \times N_2} \to \bR^m$ satisfies a rank-restricted isometry property \cite{RRIP},
e.g. by performing nuclear norm minimization \cite{NucNorm},
i.e., by outputting a solution of
\be
\label{NucNormMin}
\minimize{\vZ \in \bR^{N_1 \times N_2}} \;
\|\vZ\|_* := \sum_{k=1}^N \sigma_k(\vZ) 
\qquad \mbox{subject to }
\cA(\vZ) = \vy,
\ee
or by executing an iterative hard thresholding scheme \cite{GolMa,JMD},
i.e., by outputting the limit  of the sequence $(\vX_n)_{n \ge 0}$
defined by $\vX_0 = {\bf 0}$ and
\be
\label{ISVT}
\vX_{n+1} = H_r(\vX_n + \mu \, \cA^*(\vy - \cA \vX_n)).
\ee

Above, and in the rest of the article,\vspace{-5mm}
\begin{itemize}
\item $\sigma_1(\vZ) \ge \cdots \ge \sigma_N(\vZ)$ denote the singular values of a matrix $\vZ \in \bR^{N_1 \times N_2}$;
\item the linear map $\cA^*: \bR^{m} \to \bR^{N_1 \times N_2}$ stands for the adjoint of $\cA$ ,
defined by the property that $\langle \cA^*(\vu), \vZ \rangle_F = \langle \vu, \cA(\vZ) \rangle$ for all $\vu \in \bR^m$ and all $\vZ \in \bR^{N_1 \times N_2}$;
\item $H_r$ represents the operator of best approximation by matrices of rank at most $r$,
so that $H_r(\vZ) = \sum_{\ell=1}^r \sigma_\ell(\vZ) \vu_\ell \vv_\ell^*$ if 
$\vZ =  \sum_{\ell=1}^N \sigma_\ell(\vZ) \vu_\ell \vv_\ell^*$ is the singular value decomposition of $\vZ$.
\end{itemize}

Our focus, however,
is put on measurement maps that do not necessarily satisfy the standard rank-restricted isometry property,
but rather a modified rank-restricted isometry property featuring the $\ell_1$-norm as an inner norm.
Precisely, we consider measurement maps $\cA: \bR^{N_1 \times N_2} \to \bR^m$
with universally bounded rank-restricted isometry ratio
\be
\gamma_r := \f{\beta_r}{\alpha_r} \ge 1,
\ee
where $\alpha_r,\beta_r$ denote the optimal constants $\alpha,\beta>0$ such that 
\be
\label{MRRIP}
\alpha \|\vZ\|_F \le \| \cA(\vZ) \|_1 \le \beta \| \vZ \|_F
\qquad \qquad \mbox{for all }
\vZ \in \bR^{N_1 \times N_2} \mbox{ with } {\rm rank}(\vZ) \le r.
\ee
For instance,
it was shown in \cite{CaiZhang}\footnote{Similar statements also appear in \cite{CSV} and \cite{CC} in the cases $r=1$ and $r=2$, respectively.} that $\gamma_r \le 3$ holds
with high probability in the case of measurements obtained by rank-one projections of the form\footnote{These are called rank-one projections because $\va_i^* \vZ \vb_i$ is the Frobenius inner product $\langle \vZ, \va_i \vb_i^* \rangle_F = \tr (\vb_i \va_i^* \vZ ) $ of the matrix $\vZ$ with the rank-one matrix $\va_i \vb_i^*$.}
\be
\label{GROP}
\cA(\vZ)_i = \va_i^* \vZ \vb_i,
\qquad
i = 1,\ldots,m,
\ee 
for some independent standard Gaussian vectors $\va_1,\ldots,\va_m \in \bR^{N_1}$ and $\vb_1,\ldots,\vb_m \in \bR^{N_2}$,
provided $m \gtrsim r \max\{N_1,N_2\}$.
Moreover, 
we will briefly justify in Section \ref{SecApp} that, 
in the case of measurements given by
\be \label{MM}
\cA(\vZ)_i = \sum_{k,\ell} \cA_{i,(k,\ell)} Z_{k,\ell},
\qquad
i = 1,\ldots,m,
\ee
for some independent identically distributed mean-zero subexponential random variables $\cA_{i,(k,\ell)}$ with variance $1/m$, say,
there are constants $c,C>0$ depending on the subexponential distribution such that 
$\gamma_r \le c$ holds with high probability as soon as $m \ge C r \max\{N_1,N_2\}$.

The article \cite{CaiZhang} also showed that,
under the modified rank-restricted isometry property \eqref{MRRIP},
recovery of rank-$r$ matrices $\vX \in \bR^{N_1 \times N_2}$ from $\vy = \cA(\vX) \in \bR^m$
can be achieved by the nuclear norm minimization \eqref{NucNormMin}.
Here, we highlight the question:\vspace{-7mm}
\begin{center}
{\em 
Under the modified rank-restricted isometry property,
is it possible to achieve low-rank recovery by an iterative hard thresholding algorithm akin to \eqref{ISVT}?
}
\end{center}\vspace{-7mm}
We shall answer this question in the affirmative.
The argument resembles the one developed  in \cite{FouLec} in the context of sparse recovery, 
but some subtle matrix-analytic refinements are necessary.

The iterative hard threshoding scheme 
we propose consists in outputting the limit of the sequence $(\vX_n)_{n \ge 0}$
defined by $\vX_0 = {\bf 0}$ and
\be
\label{MISVT}
\vX_{n+1} = H_s( \vX_n + \mu_n H_t( \cA^* \sgn(\vy-\cA \vX_n) ) ),
\ee
with parameters $s$, $t$, and $\mu_n$ to be given explicitly later.
\rev{If $H_t$ did not appear, 
the scheme \eqref{MISVT} would intuitively be interpreted as a (sub)gradient descent steps for the function $\vZ \mapsto \|\vy - \cA \vZ\|_1$, followed by some singular value thresholding.}
The appearance of $H_t$ is the main adjustment distinguishing the scheme \eqref{MISVT} from the iterative hard thresholding scheme proposed in~\cite{FouLec}
(it also leads to an altered stepsize $\mu_n$).
For the more classical iterative hard thresholding scheme,
a similar adjustment was already exploited in \cite{HegIndSch},
where the roles of $H_s$ and $H_t$ were to provide `tail' and `head' approximations, respectively.
This is not what justifies the appearance of $H_t$ here.\footnote{\revv{\label{FN1}
In fact, the experiments performed in Section \ref{SecNum} suggest that $H_t$ is not necessary, but we were unable to prove the recovery of $\vX$ using the algorithm `without $H_t$'.
In the vector case, one avoids $H_t$ by relying on the observation that $H_s(\vz) = H_s(\vz_T)$
if the index set $T$ contains the support of $H_s(\vz)$.
In the matrix case, one would seemingly require the counterpart statement that
$H_s(\vZ) = H_s(P_{\cT}(\vZ))$, 
$P_\cT$ denoting the orthogonal projection onto a space $\cT$ containing the space ${\rm span}\{\vu_1\vv_1^*,\ldots,\vu_s\vv_s^*\}$ obtained from the singular value decomposition $\sum_{\ell=1}^N \sigma_\ell(\vZ) \vu_\ell \vv_\ell^*$ of $\vZ$.
Such a statement is unfortunately invalid --- take e.g. 
$\vZ = \bbmx 3 & 0 \\ 0 & 3c \ebmx$ with $|c|<1$,
so that $H_1(\vZ) = \bbmx 3 & 0 \\ 0 & 0 \ebmx$,
and take $\cT = {\rm span}\left\{ \bbmx 1 & 0\\ 0 & 0 \ebmx,
\bbmx 1 & 1 \\ 1 & 1 \ebmx \right\}$,
so that $P_\cT(\vZ) = \bbmx 3 & c\\ c & c \ebmx$
and $H_1(P_\cT(\vZ)) \not= \bbmx 3 & 0 \\ 0 & 0 \ebmx$
(otherwise ${\rm rank}\left( \bbmx 0 & c \\ c & c \ebmx \right)=1$).}}

\rev{Let us now state our most representative result.
To reiterate, it provides our newly proposed algorithm with theoretical guarantees for the success of recovery of low-rank matrices acquired via rank-one projections.
Other iterative thresholding algorithms do not possess such guarantees, 
since their analysis exploits the standard rank-restricted isometry property
and hence does not apply to rank-one projections.
We also emphasize that the result is valid beyond rank-one projections, since it only relies on the modified rank-restricted isometry property \eqref{MRRIP}.}

\bthm \label{ThmSimple}
Let $\cA: \bR^{N_1 \times N_2} \to \bR^m$ be a measurement map satisfying the modified rank-restricted isometry property \eqref{MRRIP} of order $c_0 r$ with ratio $\gamma_{c_0 r} \le \gamma$.
Then any matrix $\vX \in \bR^{N_1 \times N_2}$ of rank at most~$r$ acquired by $\vy = \cA (\vX)$
can be exactly recovered as the limit of the sequence $(\vX_n)_{n \ge 0}$ produced by \eqref{MISVT} with parameters $s = c_1 r$ , $t = c_2 r$,
and $\mu_n = \|\vy  - \cA \vX_n\|_1 / \|  H_t(\cA^* \sgn(\vy-\cA \vX_n) \|_F^2$.
The constants $c_0,c_1,c_2 >0$ depend only on $\gamma$.
\ethm

Theorem \ref{ThmSimple} is a simplified version of our main result (Theorem \ref{ThmGal}),
which covers the more realistic situation where the matrices $\vX \in \bR^{N_1 \times N_2}$ to be recovered are not exactly low-rank and where the measurements are not perfectly accurate,
i.e., $\vy = \cA \vX + \ve$ with a nonzero vector $\ve \in \bR^m$. 
The full result will be stated and proved in Section \ref{SecMain}.
In preparation, we collect in Section \ref{SecMA} the matrix-analytic tools that are essential to our argument.
In Section \ref{SecNum}, we present some modest numerical experiments demonstrating some strong points of our iterative hard thresholding algorithm.
In Section \ref{SecApp}, we close the article with a sketched justification of the modified 
rank-restricted isometry property for subexponential measurement maps.

\section{Matrix-analytic ingredients}\label{SecMA}

The proof of Theorem \ref{ThmSimple} and of its generalization (Theorem \ref{ThmGal}) is based on the two technical side results given as Propositions \ref{PropShenLi} and \ref{PropFObs} below.
The first side result provides a (loose\footnote{The expression of $\eta(\kappa)$ in \eqref{ShenLi} is certainly not optimal, but it suffices for our purpose.}) quantitative confirmation to the intuitive fact that $\|\vX - H_s(\vZ)\|_F$ approaches $\|\vX - \vZ\|_F$ as $s$ grows.
The (tight) vector version of this assertion, due to \cite{ShenLi},
was already harnessed by \cite{FouLec} in the context of sparse recovery when the standard restricted isometry property is absent.

\bprop 
\label{PropShenLi}
Given matrices $\vX,\vZ \in \bR^{N_1 \times N_2}$ and integers $r \le s$ with $s+2r \le \min\{ N_1, N_2\}$,
if ${\rm  rank}(\vX) \le r$, then
\be \label{ShenLi}
\|\vX - H_s(\vZ)\|_F \le \eta \left( \f{s}{r} \right) \|\vX - \vZ\|_F,
\qquad \mbox{where} \quad \eta(\kappa) := 1+\sqrt{\f{8}{\kappa + 1}} \underset{\kappa \to \infty}{\tto} 1. 
\ee
\eprop

The second side result says that,
under the modified rank-restricted isometry property \eqref{MRRIP},
a low-rank matrix $\vX$ is well approximated by forming $\cA^*$ applied to the sign vector of $\cA(\vX)$ and then truncating its singular value decomposition.
A closely connected result can found in \cite[Lemma~3]{FouLyn}\footnote{\cite[Lemma~3]{FouLyn} is not directly exploitable due to the obstruction described in footnote \ref{FN1}.}. 

\bprop
\label{PropFObs}
Given a matrix $\vX \in \bR^{N_1 \times N_2}$,
integers $r \le s$ with $s+r \le \min\{ N_1, N_2\}$,
and a vector $\ve \in \bR^m$,
if ${\rm  rank}(\vX) \le r$
and if the modified rank-restricted isometry property \eqref{MRRIP} holds with ratio $\gamma_{s+r} = \beta_{s+r}/\alpha_{s+r}$, then
\be
\label{BdFobs}
\|\vX - \mu H_s(\cA^*\sgn(\cA \vX + \ve)) \|_F
\le  \left( 1 - \f{1}{2\gamma_{s+r}^2} +  \tau \sqrt{\f{r}{s+r}} \right) \|\vX\|_F
+ \f{6 \tau^2}{\beta_{s+r}} \|\ve\|_1,
\ee
provided the stepsize $\mu$ satisfies, for some $\tau \ge 1$,
\be
\label{CondBdFObs}
\mu = \f{\|\cA \vX + \ve\|_1}{\nu}
\qquad \mbox{with }
\f{\beta_{s+r}^2}{\tau} \le \nu \le \beta_{s+r}^2
\mbox{ and } \nu  \ge \|H_s(\cA^*\sgn(\cA \vX + \ve))\|_F^2.
\ee
\eprop

Propositions \ref{PropShenLi} and \ref{PropFObs} both rely on a purely matrix-analytic observation stated below.

\blem
\label{LemKey}
Given matrices $\vA,\vB \in \bR^{N_1 \times N_2}$
and integers $i,j,k$ with $i \le k$ and $k+j \le \min\{N_1,N_2\}$,
if ${\rm rank}(\vA) \le j$,
then
\be
\left| \langle \vA, \vB - H_k(\vB)  \rangle_F \right|
\le
\sqrt{\f{j}{k+j-i}} \, \|\vA\|_F \, \| H_{k+j}(\vB) - H_i(\vB)  \|_F.
\ee
\elem

\bpf
We apply Von Neumman's trace inequality, 
take into account the facts that $\sigma_\ell(\vA) = 0$ for $\ell > j$
and that $\sigma_\ell(\vB - H_k(\vB)) = \sigma_{k+\ell}(\vB)$,
and finally use Cauchy--Schwarz inequality to write
\begin{align}
\label{FObs1}
\left| \langle \vA, \vB - H_k(\vB) \rangle_F \right| & 
\le \sum_{\ell \ge 1} \sigma_\ell(\vA) \sigma_\ell(\vB - H_k(\vB))\\
\nonumber
& = \sum_{\ell=1}^j \sigma_\ell(\vA) \sigma_{k+\ell}(\vB)\\
\nonumber
& \le \left[ \sum_{\ell=1}^j \sigma_\ell(\vA)^2\right]^{1/2} \left[\sum_{\ell=1}^j \sigma_{k+\ell}(\vB)^2 \right]^{1/2}.
\end{align}
For $\la_1 \ge \la_2 \ge \cdots \ge \la_N$, $N:= \min\{N_1,N_2\}$,
it is easy to verify that
\be
\sum_{\ell=1}^j \la_{k+\ell} \le \f{j}{k+j-i} \sum_{\ell=1}^{k+j-i} \la_{i+\ell}.
\ee
By applying this inequality with $\la_h = \sigma_h(\vB)^2$ and substituting into \eqref{FObs1}, we obtain
\begin{align}
\left| \langle \vA, \vB - H_k(\vB) \rangle_F \right|
& \le \|\vA\|_F \left[ \f{j}{k+j-i} \sum_{\ell=1}^{k+j-i} \sigma_{i + \ell}(\vB)^2 \right]^{1/2}\\
\nonumber
& = \sqrt{\f{j}{k+j-i}} \, \|\vA\|_F \| H_{k+j}(\vB) - H_i(\vB) \|_F, 
\end{align}
which is the announced result.
\epf

We are now in a position to fully justify Propositions \ref{PropShenLi} and \ref{PropFObs}.
 
\bpf[Proof of Proposition \ref{PropShenLi}]
We start by noticing that
\begin{align}
\|\vX - H_s(\vZ) \|_F^2  &= \langle \vX - H_s(\vZ), \vX -\vZ \rangle_F + \langle \vX - H_s(\vZ), \vZ - H_s(\vZ) \rangle_F\\
\nonumber
& \le \f{1}{2} \|\vX - H_s(\vZ) \|_F^2  + \f{1}{2} \|\vX - \vZ \|_F^2 + \langle \vX - H_s(\vZ), \vZ - H_s(\vZ) \rangle_F.
\end{align}
Thus, after rearranging the terms, we have
\begin{align}
\label{ArgShenLi1}
\|\vX - H_s(\vZ) \|_F^2  
& \le \|\vX - \vZ \|_F^2 + 2 \langle \vX - H_s(\vZ), \vZ - H_s(\vZ) \rangle_F\\
\nonumber 
& = \|\vX - \vZ \|_F^2 + 2 \langle \vX - H_r(\vZ), \vZ - H_s(\vZ) \rangle_F,
\end{align}
where the last equality followed from the fact that $\langle H_\ell(\vZ), \vZ - H_s(\vZ) \rangle_F = 0$ whenever $\ell \le s$.
We apply Lemma \ref{LemKey} with $\vA = \vX- H_r(\vZ)$, $\vB = \vZ$,
$i=r$, $j=2r$, and $k=s$.
This allows us to write 
\begin{align}
\langle \vX - H_r(\vZ), \vZ - H_s(\vZ) \rangle_F 
& \le \sqrt{\f{2r}{s+r}} \|\vX-H_r(\vZ)\|_F \| H_{s+2r}(\vZ) - H_r(\vZ) \|_F\\
\nonumber
& \le \sqrt{\f{2r}{s+r}} \|\vX-H_r(\vZ)\|_F \|\vZ - H_r(\vZ) \|_F.
\end{align}
In view of $\|\vZ - H_r(\vZ) \|_F \le \|\vZ - \vX \|_F$ (recall that $H_r(\vZ)$ is a best approximant of rank at most $r$ to $\vZ$)
and of $\|\vX-H_r(\vZ)\|_F  \le \|\vZ - H_r(\vZ) \|_F + \|\vZ - \vX \|_F \le 2 \|\vZ - \vX \|_F $, we now arrive at
\be
\label{ArgShenLi2}
\langle \vX - H_r(\vZ), \vZ - H_s(\vZ) \rangle_F 
\le 2 \sqrt{\f{2r}{s+r}}  \|\vZ - \vX \|_F^2. 
\ee
Substituting \eqref{ArgShenLi2} into \eqref{ArgShenLi1} leads to
\be
\|\vX - H_s(\vZ) \|_F^2   \le \left( 1 + 4 \sqrt{\f{2r}{s+r}}  \right)\|\vZ - \vX \|_F^2
\le \left( 1 + 2 \sqrt{\f{2r}{s+r}} \right)^2 \|\vX - \vZ \|_F^2,
\ee
which yields the desired result after taking the square root on both sides. 
\epf

\bpf[Proof of Proposition \ref{PropFObs}]
For convenience, we denote $\vw := \sgn(\cA \vX + \ve)$ throughout the proof.
Before addressing the bound \eqref{BdFobs} on $\|\vX - \mu H_s(\cA^* \vw) \|_F$,
we remark that 
$\| H_s(\cA^*\vw) \|_F^2 \le \beta_{s+r}^2$,
so that \eqref{CondBdFObs} does not feature a vacuous requirement.
In fact, a stronger inequality holds, namely 
\be
\label{ArgFObs1}
\| H_{s+r}(\cA^*\vw) \|_F \le \beta_{s+r}.
\ee
This follows from the observation that 
\begin{align}
 \| H_{s+r}(\cA^* \vw ) \|_F^2
 & = \langle H_{s+r}(\cA^* \vw), H_{s+r}(\cA^* \vw)\rangle_F
 = \langle \cA^*\vw, H_{s+r}(\cA^* \vw)\rangle_F\\
 & \nonumber
 = \langle \vw, \cA( H_{s+r}(\cA^* \vw ) )\rangle
\le \| \cA( H_{s+r}(\cA^* \vw ) ) \|_1   \\
\nonumber
& \le  \beta_{s+r} \| H_{s+r}(\cA^* \vw )\|_F.
\end{align}
Turning now to $\|\vX - \mu H_s(\cA^* \vw) \|_F$,
we expand its square to obtain
\begin{align}
\label{SquareExp}
\|\vX - \mu H_s(\cA^* \vw) \|_F^2
& = \|\vX \|_F^2 - 2 \mu \langle \vX , H_s(\cA^* \vw) \rangle_F + \mu^2 \| H_s(\cA^* \vw) \|_F^2\\
\nonumber
& \le 
\|\vX \|_F^2 - 2 \mu \langle \vX , H_s(\cA^* \vw) \rangle_F + \mu^2 \nu.
\end{align}
To deal with the inner product appearing on the right-hand side of \eqref{SquareExp}, we write
\be \label{ArgFObsAux0}
\langle \vX , H_s(\cA^* \vw ) \rangle_F
= \langle \vX , \cA^* \vw \rangle_F
- \langle \vX , \cA^*\vw  - H_s(\cA^*\vw) \rangle_F,
\ee
while noticing that
\be \label{ArgFObsAux1}
 \langle \vX , \cA^* \vw \rangle_F = \langle \cA  \vX, \sgn(\cA \vX + \ve) \rangle
 = \|\cA \vX + \ve\|_1 - \langle \ve, \sgn(\cA \vX + \ve) \rangle
 \ge \|\cA \vX + \ve \|_1 - \|\ve\|_1
\ee
and that  Lemma \ref{LemKey} applied with $\vA = \vX$, $\vB = \cA^* \vw$,
$i=0$, $j=r$, and $k=s$ yields
\be \label{ArgFObsAux2}
\left| \langle \vX , \cA^* \vw  - H_s(\cA^* \vw ) \rangle_F  \right|
 \le \sqrt{\f{r}{s+r}} \|\vX\|_F \| H_{s+r}(\cA^* \vw) \|_F.
\ee
Taking \eqref{ArgFObs1} into account
and making use of \eqref{ArgFObsAux1} and \eqref{ArgFObsAux2} in \eqref{ArgFObsAux0}, we derive
\be
\label{ArgFObs2}
\langle \vX , H_s(\cA^* \vw) \rangle_F 
\ge \|\cA \vX  + \ve \|_1 
- \left[ \|\ve\|_1 + \beta_{s+r} \sqrt{\f{r}{s+r}} \|\vX\|_F \right].
\ee
Substituting \eqref{ArgFObs2} into \eqref{SquareExp} 
and using \eqref{CondBdFObs}
now gives
\begin{align}
\|\vX - \mu H_s( \cA^* \vw ) \|_F^2
& \le \|\vX \|_F^2
- 2 \mu \|\cA \vX + \ve \|_1 + \mu^2 \nu + 2 \mu \left[ \|\ve\|_1 + \beta_{s+r} \sqrt{\f{r}{s+r}} \|\vX\|_F \right]\\
\nonumber
& = \|\vX \|_F^2 
- \f{\|\cA \vX + \ve \|_1^2}{\nu} +  2 \mu \left[ \|\ve\|_1 + \beta_{s+r} \sqrt{\f{r}{s+r}} \|\vX\|_F \right]\\
\nonumber
& \le \|\vX \|_F^2 
- \f{\|\cA \vX + \ve \|_1^2}{\beta_{s+r}^2} +  \f{2 \tau \|\cA \vX + \ve\|_1}{\beta_{s+r}^2} \left[ \|\ve\|_1 + \beta_{s+r} \sqrt{\f{r}{s+r}} \|\vX\|_F \right].
\end{align}
In view of the inequalities $\|\cA \vX + \ve\|_1^2 \ge \left| \|\cA \vX\|_1 - \|\ve\|_1 \right|^2
= \|\cA \vX\|_1^2 -2 \|\cA \vX\|_1 \|\ve\|_1 + \|\ve\|_1^2$
and  $\| \cA \vX + \ve \|_1 \le \|\cA \vX\|_1 + \|\ve\|_1$,
we arrive at
\be
\|\vX - \mu H_s( \cA^* \vw ) \|_F^2 \le a + b \|\ve\|_1 + c \|\ve\|_1^2,
\quad 
\left\{ 
\begin{matrix*}[l]
a & := \|\vX\|_F^2 - \df{\| \cA \vX  \|_1^2}{\beta_{s+r}^2} + \df{2 \tau}{\beta_{s+r}} \sqrt{\df{r}{s+r}} \|\cA \vX \|_1 \|\vX\|_F,\\
b & := \df{2(1+\tau) \|\cA \vX\|_1}{\beta_{s+r}^2} + \df{2 \tau}{\beta_{s+r}} \sqrt{\df{r}{s+r}} \|\vX\|_F,\\
c & := \df{2 \tau - 1}{\beta_{s+r}^2}.
\end{matrix*}
 \right.
\ee
The bounds $\|\cA \vX\|_1 \ge \alpha_{s+r} \|\vX\|_F$ and $\|\cA \vX\|_1 \le \beta_{s+r} \| \vX\|_F$ imply that
\be
\left\{ 
\begin{matrix*}[l]
a & \le \left( 1 - \df{1}{\gamma_{s+r}^2} +  2 \tau \sqrt{\df{r}{s+r}} \right) \|\vX\|_F^2 
\le \left( 1 - \df{1}{2\gamma_{s+r}^2} + \tau  \sqrt{\df{r}{s+r}} \right)^2 \|\vX\|_F^2 ,\\
b & \le \df{2}{\beta_{s+r}} \left( 1+\tau + \tau \sqrt{\df{r}{s+r}} \right) \|\vX\|_F
\le \df{2 \tau}{\beta_{s+r}} \left( 1+\tau + \tau \sqrt{\df{r}{s+r}} \right) \|\vX\|_F.
\end{matrix*}
 \right.
 \ee
Since we also have
\be 
c \le\f{\tau^2}{\beta_{s+r}^2}
\le \f{\tau^2 \left( 1+\tau + \tau \sqrt{\df{r}{s+r}} \right)^2}{\beta_{s+r}^2 \left( 1 - \df{1}{2\gamma_{s+r}^2} + \tau  \sqrt{\df{r}{s+r}} \right)^2},
\ee
we deduce that $\|\vX - \mu H_s ( \cA^* \vw ) \|_F^2$ is bounded above by
\be
\left[  
\left( 1 - \df{1}{2\gamma_{s+r}^2} +   \tau \sqrt{\df{r}{s+r}} \right) \|\vX\|_F
+ \f{\tau \left(1+\tau + \tau \sqrt{\df{r}{s+r}} \right)}{\beta_{s+r} \left( 1 - \df{1}{2\gamma_{s+r}^2} + \tau  \sqrt{\df{r}{s+r}} \right)} \|\ve\|_1
\right]^2.
\ee
It remains to take the square root on both sides to obtain the desired result, after noticing that $1+\tau + \tau \sqrt{r/(s+r)} \le 3 \tau$
and $1-1/(2\gamma_{s+r}^2) + \tau \sqrt{r/(s+r)} \ge 1 - 1/(2 \gamma_{s+r}^2) \ge 1/2 $.
\epf

\section{Stability and robustness of the reconstruction}\label{SecMain}

This section is devoted to the proof of Theorem \ref{ThmSimple} ---
in fact, of the more general theorem stated below,
which incorporates low-rank defect and measurement error.

\bthm \label{ThmGal}
Let $\cA: \bR^{N_1 \times N_2} \to \bR^m$ be a measurement map satisfying the modified rank-restricted isometry property \eqref{MRRIP} of order $\gamma_{t+s+r} \le \gamma$, where
\be
\label{SnT}
s = 100 \gamma^4 r
\qquad \mbox{and} \qquad
t = 800 \gamma^{12} s.
\ee
For all $\vX \in \bR^{N_1 \times N_2}$ and all $\ve \in \bR^m$,
the iterative hard thresholding algorithm \eqref{MISVT} applied to $\vy = \cA \vX + \ve$ with stepsize
\be
\label{StepSize}
\mu_n = 
\f{\|\vy- \cA \vX_n\|_1}{\max \left\{ \|H_t(\cA^* \sgn(\vy - \cA \vX_n))\|_F^2,
\df{1}{4 \gamma^2} \df{\| \cA(H_t(\cA^* \sgn(\vy - \cA \vX_n))) \|_1^2}{\|H_t(\cA^* \sgn(\vy - \cA \vX_n))\|_F^2} \right\} },
\ee
produces a sequence $(\vX_n)_{n \ge 0}$ whose cluster points $\vX^\sharp$ approximate $H_r(\vX)$ with error
\be 
\label{TheorGuarantee}
\| H_r(\vX) - \vX^\sharp \|_F
\le d \| \cA(\vX - H_r(\vX)) + \ve \|_1,
\ee
where the constant $d$ depends only on $\gamma$ and $\beta_{t+s+r}$.
\ethm

A number of comments are on order before giving the full justification of Theorem~\ref{ThmGal}.

{\bf The parameters $s$ and $t$.}
The constants involved in \eqref{SnT} are admittedly huge.
They have been chosen to look nice and to make the theory work,
but they do not reflect reality.
In practice, the numerical experiments carried our in Section \ref{SecNum} suggest that they can be taken much smaller. 

{\bf The idealized situation.}
Theorem \ref{ThmSimple} is truly a special case of Theorem \ref{ThmGal} 
where one assumes that ${\rm rank}(\vX) \le r$ and $\ve  = {\bf 0}$.
Not only does \eqref{TheorGuarantee} guarantees exact recovery in this case,
but the stepsize \eqref{StepSize} indeed reduces to $\mu_n = \|\vy - \cA \vX_n\|_1 / \| H_t(\cA^* \sgn(\vy - \cA \vX_n)) \|_F^2$
when ${\rm rank}(\vX) \le r$ and $\ve  = {\bf 0}$.
This is valid because, on the one hand,
\be
\df{1}{4 \gamma^2} \df{\| \cA(H_t(\cA^* \sgn(\vy - \cA \vX_n))) \|_1^2}{\|H_t(\cA^* \sgn(\vy - \cA \vX_n))\|_F^2}
\le \df{1}{4 \gamma^2} \beta_t^2 \le 
\df{1}{4 \gamma^2} \beta_{t+s+r}^2 \le \f{\alpha_{t+s+r}^2}{4},
\ee
and, on the other hand, $\| H_t(\cA^* \sgn(\vy - \cA \vX_n)) \|_F^2 \ge \| H_{s+r}(\cA^* \sgn(\vy - \cA \vX_n)) \|_F^2 \ge \alpha_{t+s+r}^2$, 
since
\begin{align}
\label{AsIn}
\alpha_{t+s+r} \|\vX-\vX_n\|_F
& \le \|\cA(\vX-\vX_n)\|_1  
= \langle \cA(\vX-\vX_n), \sgn(\cA(\vX-\vX_n)) \rangle\\
\nonumber
& = \langle \vX-\vX_n, \cA^*\sgn(\vy - \cA \vX_n) \rangle
\le
 \sum_{\ell = 1}^{s+r} \sigma_\ell(\vX - \vX_n)  \sigma_\ell( \cA^*\sgn(\vy - \cA \vX_n))\\
\nonumber
& \le 
\|\vX - \vX_n \|_F \|H_{s+r} (\cA^*\sgn(\vy - \cA \vX_n))\|_F.
\end{align}

{\bf The stepsize.}
The selection made in \eqref{StepSize} is certainly not the only possible one.
For instance,
choosing $\mu_n = \|\vy - \cA \vX_n \|_1 / \beta_{t+s+r}^2$ would also work.
We favored the option~\eqref{StepSize}
because it does not necessitate an estimation of $\beta_{t+s+r}$,
because it seems to perform better in practice,
and because it is close to optimal in the sense that 
the choice $\mu_n = \|\vy - \cA \vX_n\|_1 / \| H_t(\cA^* \sgn(\vy - \cA \vX_n)) \|_F^2$, made in the idealized situation, almost minimizes the expression \eqref{SquareExp} when taking \eqref{ArgFObs2} into account.
Let us also remark that there is some leeway around an admissible $\mu_n$,
namely replacing it by $\mu_n'=(1+\delta_n)\mu_n$ for a suitably small relative variation $\delta_n$ does not compromise the end result,
as it only creates a minor change in the proof (specifically, in \eqref{Comb2}).
The reader is invited to fill in the details.

{\bf Number of iterations.}
The proof of Theorem \ref{ThmGal} actually provides error estimates throughout the iterations.
For instance,
in the idealized situation where ${\rm rank}(\vX) \le r$ and $\ve = {\bf 0}$,
it reveals that there is a constant $\rho < 1$ such that
\be
\|\vX - \vX_n\|_F \le \rho^n \|\vX\|_F
\qquad \mbox{for all } n \ge 0.
\ee
Thus, for a prescribed accuracy $\eps$, we have $\|\vX - \vX_n\|_F \le \eps$ in only $n = \lceil \ln(\|\vX\|_F / \eps)/\ln(1/\rho) \rceil$ iterations.

{\bf A stopping criterion.}
Although Theorem \ref{ThmGal} states a recovery guarantee for cluster points of the sequence $(\vX_n)_{n \ge 0}$,
one can stop the iterative process as soon as the maximum in \eqref{StepSize} is achieved by the second term defining it (provided this occurs).
Precisely, if
\be
\label{CondSri}
\| \cA(H_t(\cA^* \sgn(\vy - \cA \vX_n) )) \|_1 \ge 2 \gamma \| H_t(\cA^* \sgn(\vy - \cA \vX_n) ) \|_F^2,
\ee
then the conclusion that 
\be
\| H_r(\vX) - \vX_n \|_F \le d \|\ve'\|_1,
\qquad \ve':= \cA(\vX - H_r(\vX)) + \ve
\ee
is already acquired.
Indeed, the condition \eqref{CondSri}
imposes that
\be
2 \gamma \| H_t(\cA^* \sgn(\vy - \cA \vX_n) ) \|_F^2
\le \beta_{t+s+r} \| H_t(\cA^* \sgn(\vy - \cA \vX_n) ) \|_F,
\ee
i.e., that
\be
\| H_t(\cA^* \sgn(\vy - \cA \vX_n) ) \|_F \le \f{\beta_{t+s+r}}{2 \gamma}.
\ee
We then derive (as in \eqref{AsIn}) that
\begin{align}
\label{Combb1}
\langle H_r(\vX) - \vX_n, \cA^*\sgn(\vy &- \cA \vX_n) \rangle
 \le \| H_r(\vX) - \vX_n \|_F \| H_{s+r} (\cA^* \sgn(\vy - \cA \vX_n)) \|_F\\
\nonumber
& \le \f{\beta_{t+s+r}}{2 \gamma} \| H_r(\vX) - \vX_n \|_F.
\end{align}
Moreover, keeping on mind that $\vy = \cA(H_r(\vX)) + \ve'$,
we notice that
\begin{align}
\label{Combb2}
\langle H_r(\vX) - \vX_n, \cA^*\sgn(\vy & - \cA \vX_n) \rangle
= \langle \cA( H_r(\vX) - \vX_n), \sgn(\vy - \cA \vX_n) \rangle\\
\nonumber 
& = \langle \vy - \cA \vX_n - \ve' , \sgn(\vy - \cA \vX_n)\rangle
\ge \|\vy - \cA \vX_n\|_1 - \|\ve'\|_1\\
\nonumber
& \ge \| \cA( H_r(\vX) - \vX_n ) \|_1 - 2 \|\ve'\|_1
\ge \alpha_{t+s+r} \|  H_r(\vX) - \vX_n  \|_F - 2 \|\ve'\|_1\\
\nonumber
& \ge \f{\beta_{t+s+r}}{\gamma} \|  H_r(\vX) - \vX_n  \|_F - 2 \|\ve'\|_1.
\end{align}
Combining \eqref{Combb1} and \eqref{Combb2} yields
\be
 \|  H_r(\vX) - \vX_n  \|_1 \le \f{4 \gamma}{\beta_{t+s+r}} \|\ve'\|_1.
\ee

{\bf Form of the error estimate.}
 In the realistic situation where $\vX$ is not exactly low-rank and measurement error occurs,
 one customarily sees recovery guarantees of the type
\be \label{Custom}
\|\vX - \vX^\sharp\|_F
\le \f{C}{\sqrt{r}} \min_{{\rm rank}(\vZ) \le r} \|\vX - \vZ\|_* + D \|\ve\|_1. 
\ee
Our iterative hard thresholding algorithm (with parameter $r$ replaced by $2r$)
also allows for such an estimate,
which is derived using the sort-and-split technique.
We omit the details, as the argument is quite classical and can be found in full in \cite{FouLec} for the vector case. 

{\bf The norm on the measurement error.} The appearance of $\|\ve\|_1$ in \eqref{Custom} might seem surprising,
as one usually encounters $\|\ve\|_2$ when the measurement maps satisfies the standard rank-restricted isometry property.
But there is no discrepancy here,
as this is just a matter of normalization.
Indeed, for a measurement map $\cG: \bR^{N_1 \times N_2} \to \bR^m$ generated as in \eqref{MM} with independent standard normal random variables $\cG_{i,(k,\ell)}$,
let us suppose that inaccurate measurements are given by $\vy = \cG \vX + \vu$.
On the one hand, the normalized version of $\cG$ that satisfies the standard rank-restricted isometry property is $\cA = \cG / \sqrt{m}$,
in which case the measurements take the form $\vy / \sqrt{m} = \cA \vX + \ve$ with $\ve = \vu / \sqrt{m}$,
so that $\|\ve\|_2 = \|\vu\|_2/\sqrt{m}$ appears in the error estimate.
On the other hand, the normalized version of $\cG$ that satisfies the modified rank-restricted isometry property is $\cA = \cG / m$,
in which case the measurements take the form $\vy / m = \cA \vX + \ve$ with $\ve = \vu / m$,
so that $\|\ve\|_1 = \|\vu\|_1/m$ appears in the error estimate.
The latter is in fact better than the former, by virtue of $\|\vu\|_1 \le \sqrt{m} \|\vu\|_2$ for $\vu \in \bR^m$.

\bpf[Proof of Theorem \ref{ThmGal}]
Without loss of generality,
we can assume that $\vX$ has rank at most $r$ and aim at proving that
\be \label{Obj3}
\|\vX - \vX^\sharp \|_F \le d \|\ve\|_1.
\ee
Indeed, we can in general interpret the measurements $\vy = \cA \vX + \ve$ made on $\vX$ as measurements $\vy = \cA (H_r(\vX)) + \ve'$, $\ve' := \cA(\vX - H_r(\vX)) + \ve$,
made on the matrix $H_r(\vX)$ which has rank is at most~$r$.
The exact low-rank result would then yield
\be
\|H_r(\vX) - \vX^\sharp\|_F \le d \|\ve'\|_1 
= d \| \cA(\vX - H_r(\vX)) + \ve\|_1, 
\ee
as desired.
So we assume from now on that ${\rm rank}(\vX) \le r$.
We are going to establish the existence of $\rho < 1$ such that,
for all $n \ge 0$,
\be \label{Obj1}
\|\vX - \vX_{n+1}\|_F \le \rho \|\vX - \vX_{n}\|_F + d' \|\ve\|_1,
\ee
which immediately implies that, for all $n \ge 0$,
\be \label{Obj2}
\|\vX - \vX_n\|_F \le \rho^n\|\vX\|_F + d \|\ve\|_1,
\qquad d:= \f{d'}{1-\rho}.
\ee
In particular, \eqref{Obj3} is obtained from \eqref{Obj2} by letting $n \to \infty$.
It remains to prove \eqref{Obj1} for a given integer $n \ge 0$.
Since $\vX_{n+1} = H_s( \vX_n + \mu_n H_t( \cA^* \sgn(\vy-\cA \vX_n) ) )$,
Proposition \ref{PropShenLi} implies
\be \label{Comb1}
\|\vX - \vX_{n+1}\|_F \le
\rho'
\|\vX - (\vX_n + \mu_n H_t( \cA^* \sgn(\vy-\cA \vX_n) )) \|_F,
\qquad \rho':=1+\sqrt{\f{8r}{s+r}} >1.
\ee
Next, we notice that the stepsize $\mu_n$ takes the form
\be
\mu_n = \f{\|\cA(\vX-\vX_n)+\ve\|_1}{\nu_n},
\ee
where the denominator 
\be
\nu_n := \max\left\{ \|H_t(\cA^* \sgn(\cA(\vX - \vX_n) + \ve))\|_F^2,
\f{1}{4 \gamma^2}\f{\| \cA (H_t(\cA^* \sgn(\cA(\vX - \vX_n) + \ve)))  \|_1^2 }{ \|H_t(\cA^* \sgn(\cA(\vX - \vX_n) + \ve))\|_F^2}
\right\} 
\ee
clearly satisfies
\be
\nu_n \ge \|H_t(\cA^* \sgn(\cA(\vX - \vX_n) + \ve))\|_F^2, 
\ee
as well as 
\be
\f{\beta_{t+s+r}^2}{4 \gamma^4} \le \nu_n \le \beta_{t+s+r}^2.
\ee
The latter inequalities follow from
$\|H_t(\cA^* \sgn(\cA(\vX - \vX_n) + \ve))\|_F^2 \le \beta_t^2 \le \beta_{t+s+r}^2$ 
(obtained in the same way as \eqref{ArgFObs1})
and from
\be
\df{1}{4 \gamma^2}\f{\| \cA (H_t(\cA^* \sgn(\cA(\vX - \vX_n) + \ve)))  \|_1^2 }{ \|H_t(\cA^* \sgn(\cA(\vX - \vX_n) + \ve))\|_F^2}
\left\{
\bmx
\le  \df{1}{4 \gamma^2} \beta_{t+s+r}^2 \le  \beta_{t+s+r}^2,\\ \\
\ge  \df{1}{4 \gamma^2} \alpha_{t+s+r}^2 \ge \df{\beta_{t+s+r}^2}{4 \gamma^4}.
\emx
\right.
\ee
Therefore, 
according to Proposition \ref{PropFObs}, we obtain
\begin{align} \label{Comb2}
\|\vX - (\vX_n & + \mu_n H_t( \cA^* \sgn(\vy-\cA \vX_n) )) \|_F
= \|(\vX - \vX_n) - \mu_n H_t( \cA^* \sgn(\cA(\vX- \vX_n) + \ve) )) \|_F\\
\nonumber
& \le \rho'' \|\vX-\vX_n\|_F
+ \f{96 \gamma^8}{\beta_{t+s+r}} \|\ve\|_1,
\qquad 
\rho'' := 1 - \f{1}{2\gamma^2} +  4 \gamma^4 \sqrt{\f{s+r}{t+s+r}} < 1,
\end{align} 
with $4 \gamma^4 \sqrt{(s+r)/(t+s+r)} < 1/(2\gamma^2)$ being justified by our choice of $s$ and $t$.
In fact, our choice of $s$~and~$t$ guarantees that
\begin{align}
\rho := \rho' \rho''
& \le \left(1+\sqrt{\f{8r}{s}} \right) \left( 1 - \f{1}{2\gamma^2} + 4 \gamma^4 \sqrt{\f{2s}{t}} \right)
\le 1 - \f{1}{2\gamma^2} + 4 \gamma^4 \sqrt{\f{2s}{t}} +\sqrt{\f{8r}{s}}\\
\nonumber
& \le 1 - \f{1}{2\gamma^2} +  \f{1}{5\gamma^2} + \f{\sqrt{2}}{5\gamma^2}
= 1- \f{3-2\sqrt{2}}{10 \gamma^2}
< 1.
\end{align}
Combining \eqref{Comb1} and \eqref{Comb2}, we finally derive that
\be
\|\vX - \vX_{n+1}\|_F \le \rho \|\vX - \vX_{n}\|_F + d' \|\ve\|_1,
\qquad d':= \f{96 \rho' \gamma^8}{\beta_{t+s+r}}.
\ee
This establishes \eqref{Obj1} and completes the proof.
\epf

\section{Numerical validation}\label{SecNum}

\revv{The perspective adopted in this article is mostly theoretical,
with a declared goal consisting in uncovering  iterative thresholding algorithms that provide alternatives to nuclear norm minimization when the standard rank-restricted isometry property does not hold.
The modest experiments carried out with our proposed algorithm are encouraging,
although we shall refrain from making too optimistic claims about its performance.}
In what follows, we only focus on measurements obtained by Gaussian rank-one projections of type \eqref{GROP}.
Note that in this case our algorithm starts with $\vX_0 = {\bf 0}$
and iterates the scheme
\be
\label{MIHT-GRk1}
\vX_{n+1} = H_s\left(  
\vX_n + \mu_n H_t \left( \sum_{i=1}^m \sgn(y_i - \va_i^* \vX_n \vb_i) \vb_i \va_i^* \right)
\right).
\ee
All the experiments performed below can be reproduced by downloading the {\sc matlab} file linked to this article from the first author's webpage.

\begin{figure}[h]
\center
\includegraphics[width=0.55\textwidth]{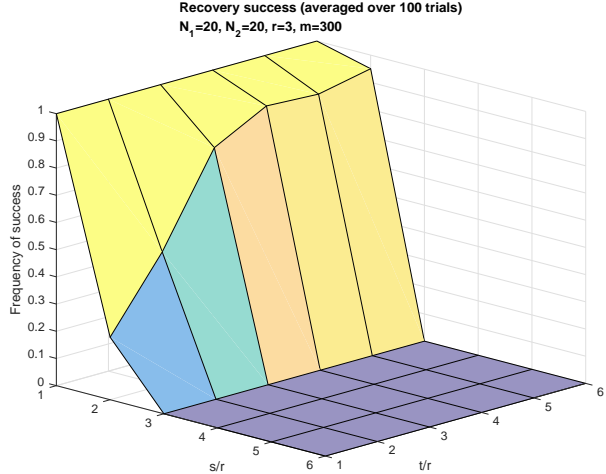}
\caption{Influence of $s$ and $t$ on the success of algorithm \eqref{MIHT-GRk1}.}
\label{FigExp1}
\end{figure}

{\bf Influence of $s$ and $t$.}
In order to derive theoretical guarantees about the success of the algorithm,
the parameters $s$ and $t$ appearing in \eqref{MIHT-GRk1} can be chosen according to \eqref{SnT}.
This of course does not mean that other choices are unsuitable.
Figure \ref{FigExp1} reports on the effect of varying these parameters.
We observe that successful recoveries occur more frequently when $s=r$,
which makes intuitive sense since the rank-$s$ matrix $\vX_n$ is supposed to approximate the rank-$r$ matrix $\vX$.
\revv{We also observe that successful recoveries occur more frequently when $t$ increases.
The largest value $t$ can take is $N=\min\{N_1,N_2\}$,
in which case the hard thresholding operator does not do anything.
The default value of the parameters are therefore $s=r$ and $t=N$,
meaning that $H_t$ is absent and that one singular value decomposition is saved per iteration.}

\begin{figure}[h]
\center
\subfigure{
\includegraphics[width=0.48\textwidth]{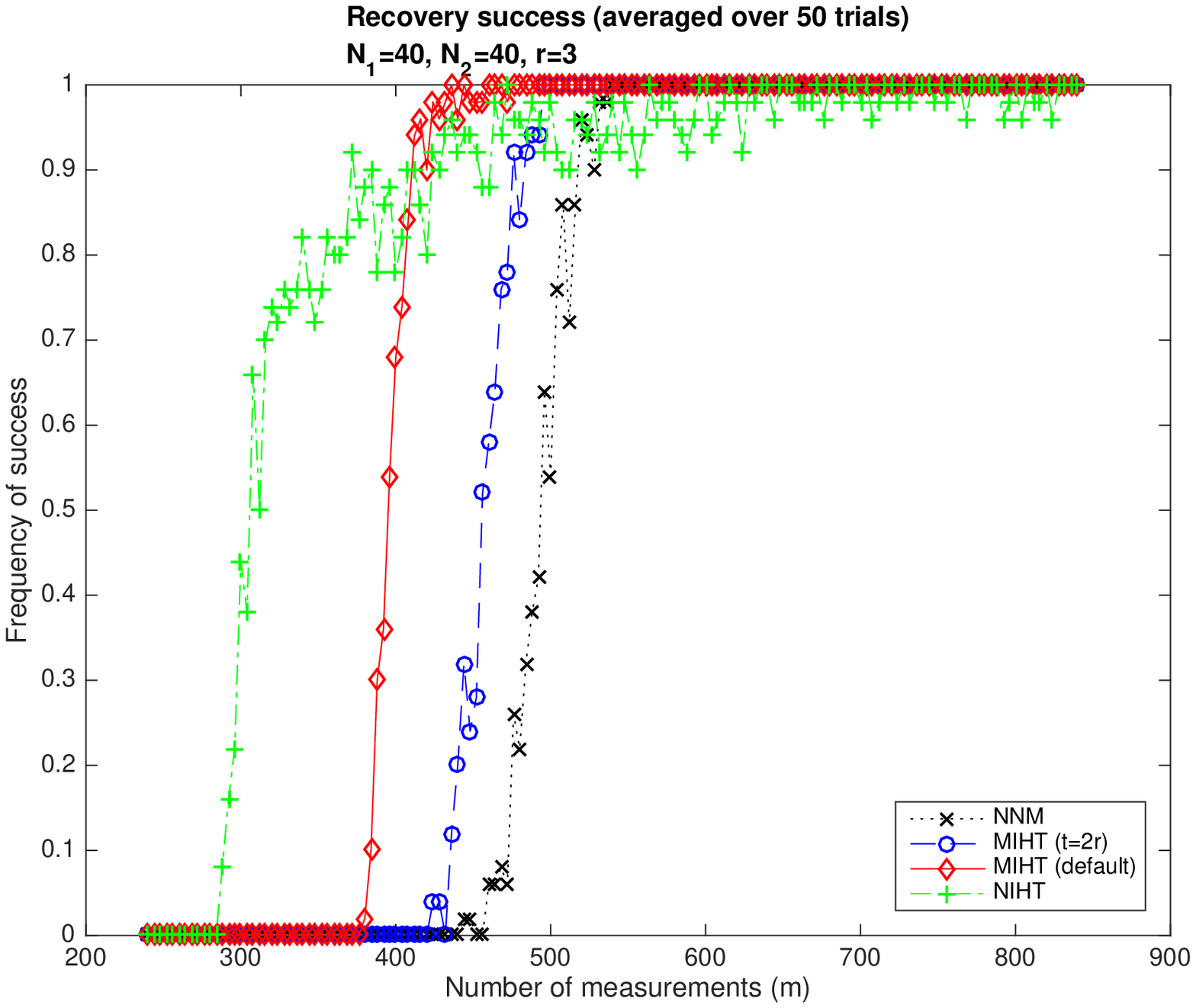}
}
\subfigure{
\includegraphics[width=0.48\textwidth]{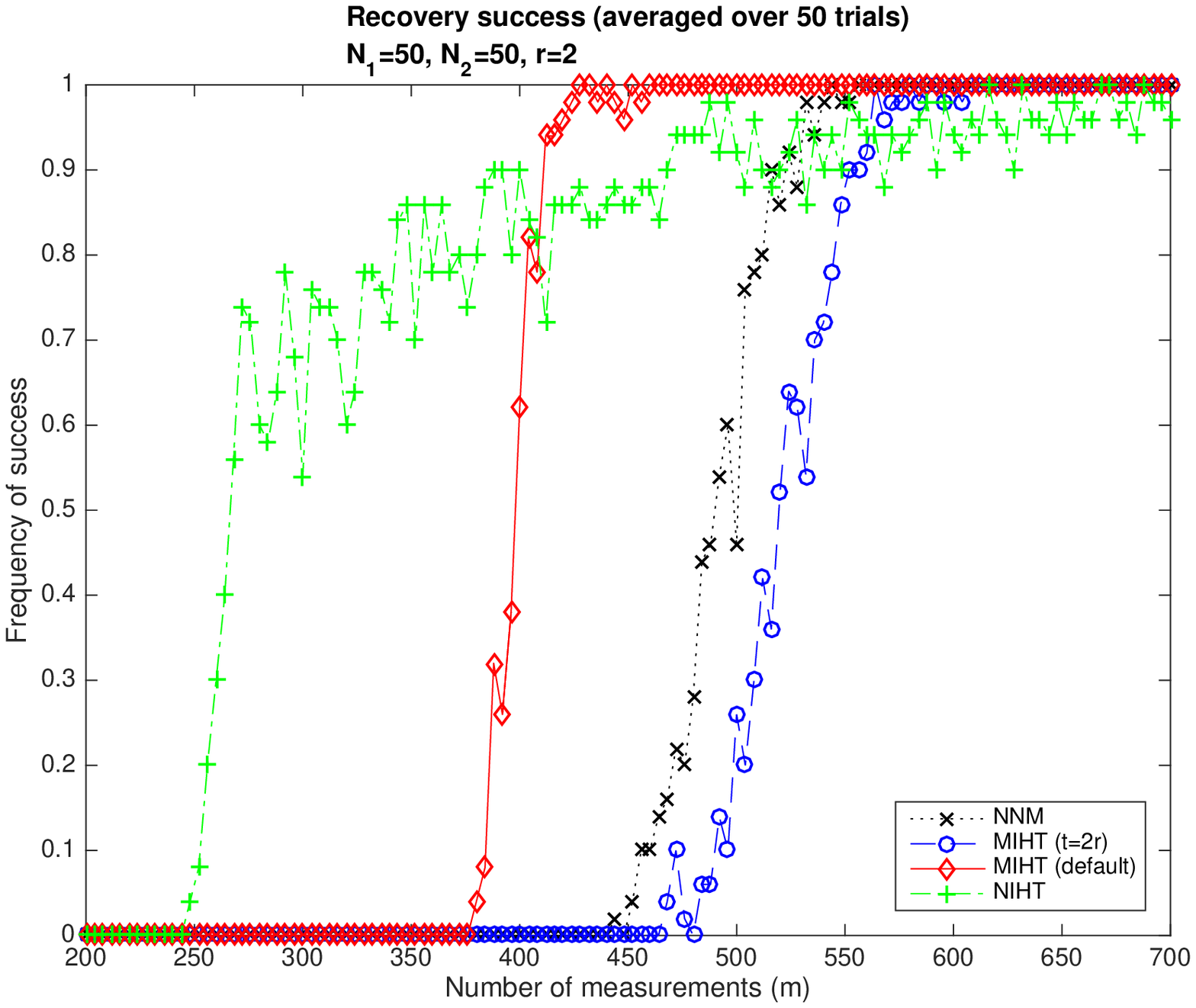}
}
\caption{\rev{Frequency of successful recoveries as a function of the number of Gaussian rank-one measurements
for nuclear norm minimization,
normalized iterative hard thresholding,
and modified iterative hard thresholding with default parameters and with parameters $s=r$ and $t=2r$.}}
\label{FigExp2}
\end{figure}

{\bf Performance.}
We now compare the low-rank recovery capability of the algorithm \eqref{MIHT-GRk1} with two other algorithms,
namely with nuclear norm minimization (NNM) 
and with the normalized iterative hard thresholding (NIHT) algorithm proposed in \cite{TanWei}.
In fact, we test two versions of our modified iterative hard thresholding (MIHT) algorithm:
one with the default parameters as discussed above,
and one with parameters $s=r$ and $t=2r$
--- this is another natural choice,
as $H_t(\cA^* \sgn(\vy - \cA \vX_n) )$ is supposed to approximate the rank-$2r$ matrix $\vX-\vX_n$.
Judging from Figure \ref{FigExp2}, NNM is seemingly (and surprisingly) outperformed by our default MIHT algorithm.
The comparison with NIHT is less straightforward.
With few measurements available, NIHT does best,
but when the number of measurements is large enough for the recoveries by NNM and MIHT to be certain,
\revv{recovery by NIHT does not occur consistently.
We interpret this phenomenon as a reflection of the theory.
Indeed, Gaussian rank-one projections of type \eqref{GROP} do not obey the standard rank-restricted isometry property,
which would have guaranteed the success of NIHT (see \cite{TanWei}),
but they obey the modified rank-restricted isometry property \eqref{MRRIP},
which does guarantee the success of NNM (see \cite{CaiZhang})
and of MIHT (this article).}

\begin{figure}[h]
\center
\subfigure{
\includegraphics[width=0.48\textwidth]{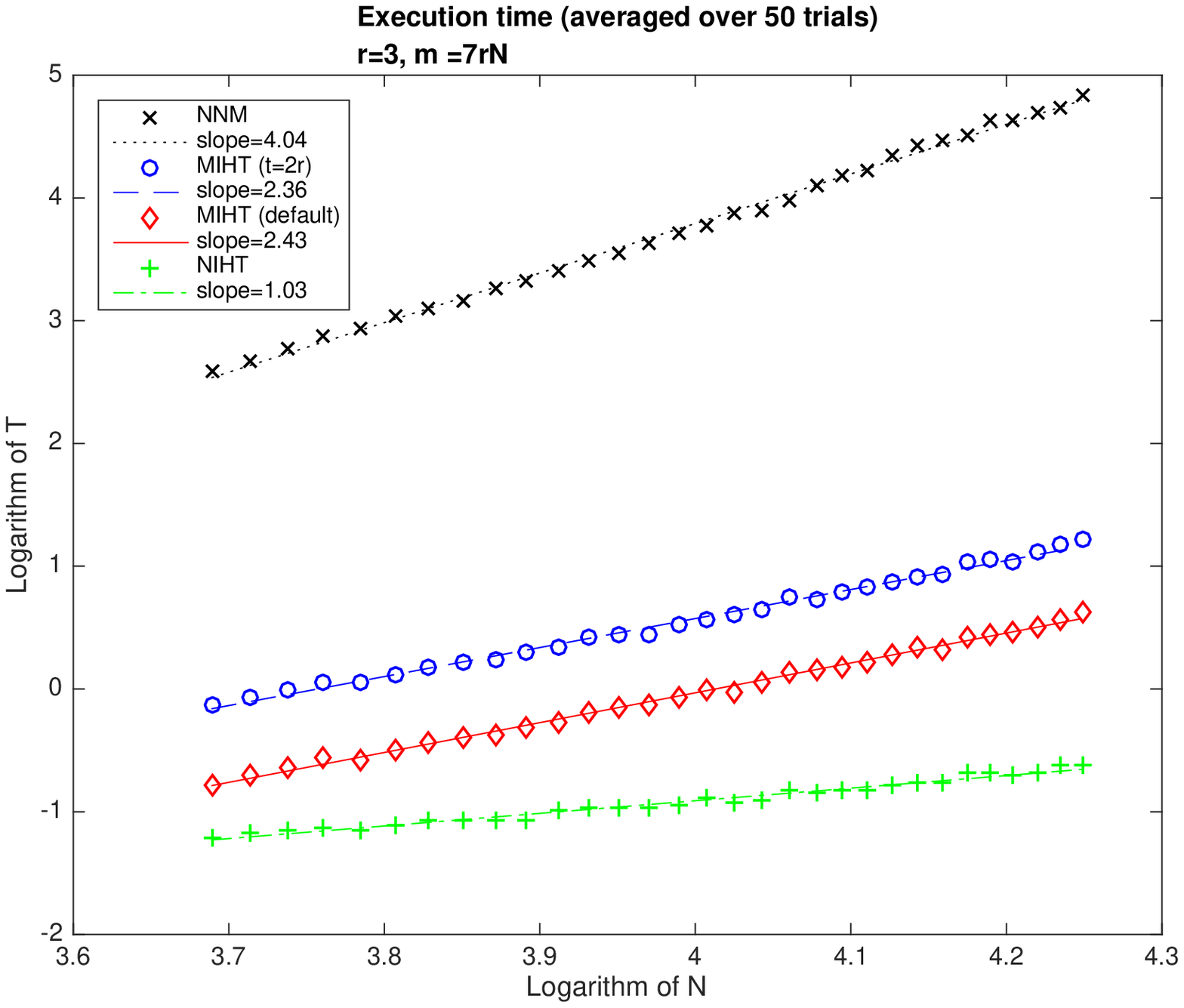}
}
\subfigure{
\includegraphics[width=0.48\textwidth]{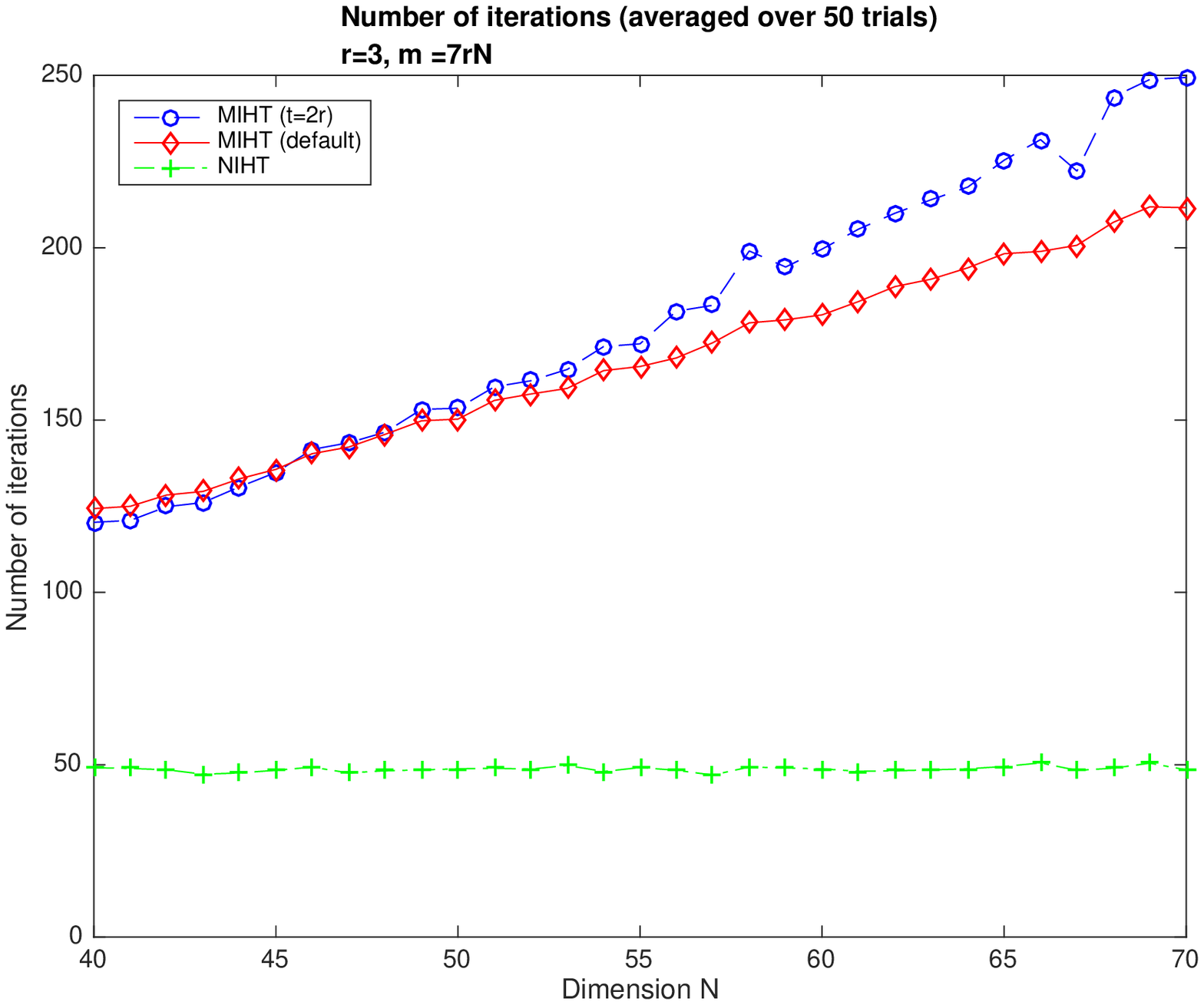}
}
\caption{Execution time and number of iterations as a function of the dimension $N=N_1=N_2$.}
\label{FigExp3}
\end{figure}

{\bf Speed.}
The next experiment compares the four algorithms previously considered in terms of execution time and number of iterations (when applicable).
Figure \ref{FigExp3} shows, as expected,
that iterative algorithms are much faster than nuclear norm minimization:
in this experiment (which only takes successful recoveries into account),
the execution time is $\Theta(N^\omega)$
with $\omega \approx 4$ for NNM,
$\omega \approx 1$ for NIHT,
and $\omega \approx 2.4$ for MIHT.
The execution times for the two versions of MIHT essentially differ by a factor of two,
due to the fact that the default version performs only one singular value decomposition per iteration instead of two.
The number of iterations is roughly similar for the two versions.
It is much lower for NIHT,
which explains its superiority in terms of overall execution time.

\begin{figure}[h]
\center
\includegraphics[width=0.55\textwidth]{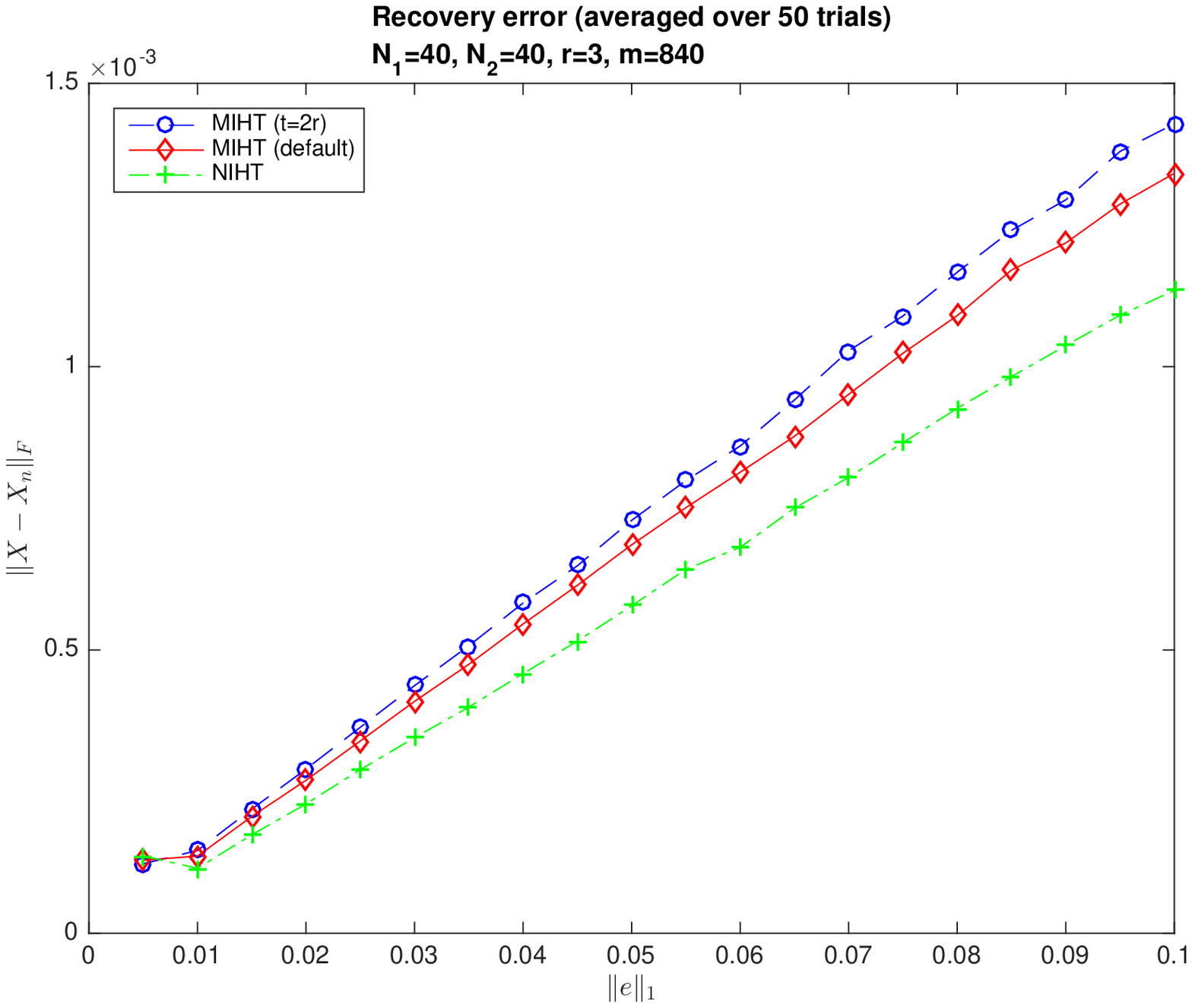}
\caption{Frobenius norm of the recovery error as a function of the $\ell_1$-norm of the measurement error.}
\label{FigExp4}
\end{figure}

{\bf Robustness.} The final experiment is designed to examine the effect on the recovery of errors $\ve \in \bR^m$ in the measurements $\vy =\cA(\vX) + \ve$ made on low-rank matrices $\vX$.
In particular, we wish to certify the statement  \eqref{TheorGuarantee} that the recovery error, measured in Frobenius norm,
is at most proportional to the measurement error, measured in $\ell_1$-norm.
The results displayed in Figure \ref{FigExp4} confirm the validity of this statement for our algorithm \eqref{MIHT-GRk1} with default parameters and with parameters $s=r$ and $t=2r$,
as well as for the normalized iterative hard thresholding algorithm.

\section{Explanation of the modified rank-restricted isometry property}\label{SecApp}

In this appendix, \revv{we point out that the modified rank-restricted isometry property \eqref{MRRIP} is not only valid for rank-one Gaussian measurements,
as already established in \cite{CaiZhang},
but also for subexponential measurements of type \eqref{MM}.}
The justification essentially follows the arguments put forward in \cite{FouLai}.
First, 
fixing a matrix $\vZ \in \bR^{N_1 \times N_2}$,
viewed as its vectorization $\vz = {\rm vec}(\vZ) \in \bR^{N_1N_2}$,
\cite[Theorem~3.1]{FouLai} guarantees a concentration inequality of the form
\be
\label{CI}
{\rm Pr}\left( \left| \|\cA (\vZ)\|_1 - \sslash \vZ \sslash  \right| > \eps \sslash \vZ \sslash \right) \le 2 \exp(- \kappa \eps^2 m).
\ee
The constant $\kappa$ depends on the subexponential distribution,
and so does the slanted norm,
although it is comparable to the Frobenius norm in the sense that (see \cite[Proposition~2.3]{FouLai})
\be
\label{SlF}
c \|\vZ\|_F \le \sslash \vZ \sslash \le C \|\vZ\|_F
\qquad
\mbox{for all } \vZ \in \bR^{N_1 \times N_2}
\ee
for some constants $c,C >0$ depending on the subexponential distribution.
The next step (see the proof of \cite[Theorem~4.1]{FouLai}) is a covering argument
showing that \eqref{CI} extends to a uniform concentration inequality of the form 
 \be
\label{UCI}
{\rm Pr}\big( \big| \|\cA (\vZ)\|_1 - \sslash \vZ \sslash  \big| > \delta \sslash \vZ \sslash 
\mbox{ for some $\vZ $ with } {\rm rank}(\vZ) \le r\big) \le 2 \exp(- \kappa \bar{\delta}^2 m + \ln( \cN_r(\bar{\delta}))),
\ee
where $\cN_r(\bar{\delta})$, $\bar{\delta} := \delta/3$,
denotes the size of a $\bar{\delta}$-net for the set of rank-$r$ matrices with slanted norm at most $1$.
We now recall the key observation from \cite{RRIP} that the set of rank-$r$ matrices in $\bR^{N_1 \times N_2}$ with norm at most $1$ has a covering number $\cN_r(\bar{\delta}) \le (c_1/\bar{\delta})^{c_2 r N_{\max}}$,
$N_{\max} := \max \{ N_1,N_2 \}$,
which is valid for the Frobenius norm and, in view of \eqref{SlF},  for the slated norm, too.
This implies a modified rank-restricted isometry property of the form
\be
\label{UCI2}
(1-\delta) \sslash \vZ \sslash \le \| \cA (\vZ) \|_1 \le (1+\delta) \sslash \vZ \sslash
\qquad \quad \mbox{whenever } {\rm rank}(\vZ) \le r,
\ee
which occurs with failure probability at most
\be
2 \exp(- \kappa \bar{\delta}^2 m + \ln(c_1/\bar{\delta})c_2 r N_{\max}  )
\le 2 \exp(- \kappa \bar{\delta}^2 m / 2 )
\le 2 \exp(- \kappa' \delta^2 m),
\ee
provided $\ln(c_1/\bar{\delta})c_2 r N_{\max} \le  \kappa' \bar{\delta}^2 m / 2 $,
i.e., $m \ge C(\delta) r N_{\max}$.
Finally, \eqref{UCI2} turns into \eqref{MRRIP} by virtue of \eqref{SlF}.

\end{document}